\documentstyle{amlts}
\begin{document}
\annalsline{155}{2002}
\received{May 9, 2001}
\startingpage{929}
\def\bye{\end{document}}
 \font\tenrm=cmr10

\catcode`\@=11
\font\twelvemsb=msbm10 scaled 1100
\font\tenmsb=msbm10
\font\ninemsb=msbm10 scaled 800
\newfam\msbfam
\textfont\msbfam=\twelvemsb  \scriptfont\msbfam=\ninemsb
  \scriptscriptfont\msbfam=\ninemsb
\def\msb@{\hexnumber@\msbfam}
\def\Bbb{\relax\ifmmode\let\next\Bbb@\else
 \def\next{\errmessage{Use \string\Bbb\space only in math
mode}}\fi\next}
\def\Bbb@#1{{\Bbb@@{#1}}}
\def\Bbb@@#1{\fam\msbfam#1}
\catcode`\@=12

 \catcode`\@=11
\font\twelveeuf=eufm10 scaled 1100
\font\teneuf=eufm10
\font\nineeuf=eufm7 scaled 1100
\newfam\euffam
\textfont\euffam=\twelveeuf  \scriptfont\euffam=\teneuf
  \scriptscriptfont\euffam=\nineeuf
\def\euf@{\hexnumber@\euffam}
\def\frak{\relax\ifmmode\let\next\frak@\else
 \def\next{\errmessage{Use \string\frak\space only in math
mode}}\fi\next}
\def\frak@#1{{\frak@@{#1}}}
\def\frak@@#1{\fam\euffam#1}
\catcode`\@=12


\newcommand \si {\Sigma_m}


\newcommand {\R} {{\Bbb R}}

\newcommand{\peta} {{\bf P}_{\eta}}

\newcommand {\sigp} {\sigma_P}

\newcommand {\sigca} {\Sigma_{\cal{A}}}

\newcommand {\siax} {\Sigma_{\cal{A}}\times X}

\newcommand {\B} {{\cal B}}
\newcommand {\G} {{\cal G}}


\newcommand \F  {{\cal F}}
\newcommand \A  {{\cal A}}
\newcommand \I  {{\cal I}}
\newcommand \yy {Y^{\Bbb Z}}
\newcommand \om {\Omega_m}
\newcommand \p {\varphi}
\newcommand \pb {\overline\varphi}
\newcommand \lox {L_1(\Omega_m\times X, \mu\times\nu)}
\newcommand \ly {L_1(Y, \eta)}
\newcommand \lx {L_1(X,\nu)}
\newcommand \lz {L_1(Z,\mu)}
\newcommand \ltwo {L_2(X,\nu)}
\newcommand \xn {X^{{\Bbb N}}}
\newcommand \ppk {E(\p_k|{\cal F}_{\infty})}
\newcommand \ppsk {E(\psi_k|{\cal F}_{\infty})}
\newcommand \Sia {\Sigma_{{\cal A}}}
\newcommand \sia {\sigma_{{\cal A}}}
\newcommand {\sn} {s_n^{\Pi}}
\newcommand {\qz} {{\bf Q}_{\mu}}
\newcommand {\pz} {{\bf P}_{\eta}}
\newcommand {\bom} {{\bf \omega}}
\newcommand {\sip} {{\sigma}_{P}}

\newcommand {\y} {{\bf y}}
\newcommand {\z} {{\bf z}}
\newcommand {\yz} {Y^{\Bbb Z}}
\newcommand {\zz} {Z^{\Bbb Z}}

\title{Convergence of spherical averages\\ for actions of free groups}
 \shorttitle{Convergence of spherical averages} 
\author {Alexander I. Bufetov}

 \institutions{Princeton University, Princeton, NJ and Independent University of Moscow, Moscow, Russia\\
{\eightpoint {\it E-mail address\/}: bufetov@math.princeton.edu}}

\section{Introduction}
\label{intro}

Let $(X, \nu)$ be a probability space and suppose a free group $F_m$ 
with $m$ generators acts on $(X,\nu)$ by measure-preserving 
transformations.
Let $\{a_1,\dots, a_m\}$ be 
 a set of free generators for 
$F_m$ and let
  $T_1,\dots, T_m: X\to X$ be transformations corresponding to the generators.
Write $T_{-i}=T_i^{-1}$ for $i=1,\dots,m$, and set 
$\A=\{-m,\dots,-1,1,\dots,m\}$. We also have the action $F_m$ on 
$L_1(X,\nu)$, defined by $T_g\p=\p\circ T_{g^{-1}}$, $g\in F_m$.

Consider the set $W_{\A}$ of all finite words over the alphabet  $\A$:
$$
W_{\A}=\{w=w_1w_2\dots w_n|\ w_i\in\A\} .
$$

Denote by $|w|$
the length of the word $w$ and for any positive integer $n$, let 
$W_{\A}(n)=\{w\in W_{\A}, |w|=n\}$.

For each $w\in W_{\A}$, $w=w_1\dots w_n$,
define a transformation
\begin{equation}
\label{fw}
T_w=T_{w_1}T_{w_2}\dots T_{w_n}.
\end{equation}

Let $\Pi$ be a stochastic $2m\times 2m$ matrix, 
whose rows and columns are indexed by elements of $\A$, that is, 
$\Pi=(p_{ij}), i,j\in\A$.
 Assume that $\Pi$ has a unique stationary distribution 
$(p_{-m},\dots,p_{-1}, p_1, \dots, p_m)$ and that   $p_i>0$ for all $i\in \A$.

For $w\in W_{\A}$, $w=w_1\dots w_n$, denote 
$$
p(w)=p_{w_nw_{n-1}}p_{w_{n-1}w_{n-2}}\dots p_{w_{2}w_1},\ \pi(w)=p_{w_n}p(w).
$$

Consider the operators 
\begin{equation}
s_n^{\Pi}=\sum_{|w|=n} \pi(w)T_w.
\end{equation}

In this paper, we investigate convergence of this sequence of operators.

\demo{Definition {\rm 1}}
We shall say that the matrix  $\Pi$ {\it generates the free group} 
if $p_{ij}=0$ is equivalent to $i+j=0$. 
\enddemo

We shall need  the {\it symmetry condition}
\begin{equation}
\label{sym}
p_i=p_{-i}, \ \ p_{-i,-j}=\frac{p_jp_{ji}}{p_i}.
\end{equation}

Relation (\ref{sym})
is equivalent to saying that all operators $\sn$ are self-adjoint.
\vglue2pt
Let $F_m^2$ be the subgroup of words of even length in 
$F_m$, that is, the subgroup generated by $a_ia_j$, 
$i,j\in\{1,\dots,m\}$. 

Recall that  $L\log L(X,\nu)=\{\p\in\lx: \int_X|\p|\log^+|\p|d\nu<\infty\}$.

\specialnumber{1}\proclaim{Theorem}
\label{sphfr}
 Let $(X,\nu)$ be a Lebesgue probability space{\rm .}
Assume the matrix $\Pi$ generates the free group
and satisfies {\rm (\ref{sym}).} Then for any
 $\p\in L\log L(X,\nu)${\rm ,} the sequence $s_{2n}^{\Pi}\p$ converges 
as $n\to\infty$ both 
$\nu$\/{\rm -}\/almost 
everywhere and in $\lx$ to an $F_m^2$\/{\rm -}\/invariant function{\rm .}
\endproclaim

{\it Remark}.
The sequence $s_{2n+1}^{\Pi}\p$ also converges. 
The sequence $s_{n}^{\Pi}$ need not converge, however, because
the action of $F_m$ might have {\it an eigenfunction with eigenvalue} $-1$, 
that is, a nonzero function
$\psi\in\lx$ such that $T_i\psi=-\psi$ for all $i\in\A$ (for the 
same reason, the limit in Theorem \ref{sphfr} must be  $F_m^2$-invariant 
but need not be $F_m$-invariant).
If the action does not have eigenfunctions with eigenvalue $-1$ 
then for any $\p\in L\log L(X,\nu)$ the sequence 
$s_n^{\Pi}\p$ converges as $n\to\infty$ both 
$\nu$-almost everywhere and in $\lx$ to 
an $F_m$-invariant limit. 
\vglue3pt

Averages $s_{2n}^{\Pi}$ converge under weaker assumptions       
on the matrix $\Pi$ than in Theorem \ref{sphfr}.

\vglue4pt {\it Definition {\rm 2}}. A
matrix  $\Pi$ with nonnegative entries will be
called {\it irreducible}
if for some $n>0$ all entries of the matrix $\Pi+\Pi^2+\dots \Pi^n$
are positive
(if $\Pi$ is stochastic then this is
equivalent to saying that in the corresponding
Markov chain any state is attainable from any other state).
\vglue4pt

 {\it Definition}  3. A matrix $\Pi$ with nonnegative entries
will be called {\it  strictly
irreducible} if $\Pi$ is irreducible and $\Pi\Pi^T$ is irreducible
(here $\Pi^T$
stands for the transpose of $\Pi$.)
\vglue4pt

Clearly, a matrix generating the free 
group is strictly irreducible.

\specialnumber{2}\proclaim{Theorem}
\label{sphir}
Let $(X,\nu)$ be a Lebesgue probability space and let $p>1${\rm .}
Assume the matrix $\Pi$ is strictly irreducible
and satisfies {\rm (\ref{sym}).} Then for any
 $\p\in L_p(X,\nu)${\rm ,} the sequence $s_{2n}^{\Pi}\p$ converges
as $n\to\infty$ both 
$\nu$\/{\rm -}\/almost everywhere and in $L_p$ to an $F_m^2$\/{\rm -}\/invariant function{\rm .}
\endproclaim

\vglue-9pt
\section{History}
\label{hist}
\vglue-9pt

\advance\eqcount by 3

First ergodic theorems for actions of arbitrary countable groups 
were obtained by V.  I.\  Oseledets \cite{oseled} in the following setting.

Let $\Gamma$ be a countable group that acts by measure-preserving 
transformations of a probability space $(X,\nu)$, and for $g\in\Gamma$
let $T_g$ be the corresponding transformation.
Let $\mu$ be a probability measure on $\Gamma$ satisfying the condition
$\mu(g^{-1})=\mu(g)$. Let $\mu^{(n)}$ be the $n$-th convolution of $\mu$. 
The ergodic theorem of Oseledets states that 
for $\p\in L\log L(X,\nu)$, the averages 

$$
A_{2n}\p=\sum_{g\in\Gamma}\mu^{(2n)}(g)T_g\p
$$
converge almost everywhere. The proof is based on consideration of
the self-adjoint Markov operator $Q=\sum_{g\in\Gamma}\mu(g)T_g$.

In 1969 Y. Guivarc'h \cite{guiv} (motivated by the work of 
Arnold and Krylov \cite{arkr}) considered uniform spherical averages 
on the free group; that is,
\begin{equation}
\label{unif}
s_n=\frac1{2m(2m-1)^{n-1}}\sum_{g:|g|=n}T_g
\end{equation}
and proved that for $\p\in L_2(X,\nu)$ the sequence $s_{2n}\p$ 
converges in $L_2$ to an $F_m^2$-invariant function. 

In 1986 R.\ I.\  Grigorchuk \cite{gri1} (see also \cite{gri2}) 
announced pointwise convergence for the averages 
$$
C_N=\frac1N\sum_{n=0}^{N-1}s_n.
$$

In 1994 Nevo and Stein proved: 

\proclaimtitle{Nevo and Stein \cite{nest}}
\specialnumber{3}\proclaim{Theorem}
Let $p>1${\rm .} Then for any $\p\in L_p(X,\nu)$ the sequence 
$s_{2n}\p$ converges as $n\to\infty$ 
both $\nu$\/{\rm -}\/almost everywhere and in $L_p$
to an $F_m^2$\/{\rm -}\/invariant function{\rm .} 
\endproclaim

The Nevo-Stein theorem is a particular case of Theorem \ref{sphfr};
we shall however consider it separately in Section \ref{secunif}
in order to illustrate the ideas of the proof of Theorem \ref{sphfr}.

\vglue-9pt
\section {The Markov operator}
\advance\eqcount by 4
\vglue-3pt

Recall that if $(Z,\mu)$ is a probability space then a linear operator
$Q$ on $L_1(Z,\mu)$ is called a {\it measure-preserving Markov operator}
if it preserves the cone of nonnegative functions, $L_1$-norm, 
and $L_{\infty}$-norm.

Let $p=\{p_{-m}, \dots, p_{-1}, p_1,\dots,p_m\}$ 
be the stationary distribution of the matrix $\Pi$.

Consider the space $Y=X\times\A$ with the measure $\eta=\nu\times p$ and
a Markov operator $P$ on $\ly$ given by 
\begin{equation}
\label{operator}
P\p(x,i)= \sum_{j\in\A}p_{ij}\p(T_ix,j).
\end{equation}

$P$ is a measure-preserving Markov operator on $L_1(Y,\eta)$. It was 
introduced 
by R.\ I.\  Grigorchuk \cite{gri2}, 
J.-P. Thouvenot (oral communication), 
and myself \cite {buf}.

For $n>1$ we have
\begin{equation}
\label{t}
P^n\p(x,i)=\sum_{w\in W_{\A}(n-1), j\in\A}p_{iw_{n-1}}p(w)
p_{w_1j}\p(T_wT_{i}x,j), 
\end{equation}
which implies:

\specialnumber{1}\proclaim{Proposition}
Let $\psi\in\lx${\rm .} Let $\p\in\ly$ be given by $\p(x,a)=\psi(x)${\rm .}
Then 
$$
\sn\psi=\sum_{i\in\A}p_iP^n\p(x,i).
$$
\endproclaim

To  prove Theorem \ref{sphfr}, it suffices to  prove the following:
\specialnumber{1}\proclaim{Lemma}
\label{ppow}
Suppose $\Pi$ generates the free group and 
satisfies the symmetry condition {\rm (\ref{sym}).}
Suppose the action of $F_m^2$ on $(X,\nu)$ is ergodic{\rm .}
Then for any $\p\in L\log L(Y,\eta)${\rm ,}
$$
P^n\p\to\int\limits_{Y}\p d\eta
$$
both $\eta$\/{\rm -}\/almost everywhere and in $\ly${\rm .}
\endproclaim

First we discuss ergodicity of $P$ and $P^2$.

\specialnumber{2}\proclaim{Lemma}
\label{erg}
If the action of $F_m$ on $(X,\nu)$ is ergodic and 
$\Pi$ is strictly irreducible{\rm ,} then $P$ is ergodic{\rm .}
\endproclaim

 {\it Definition} {\rm 4}.
A function $\p\in\ly$ {\it does not
depend on} $\A$ if there exists $\psi\in\lx$ such that $\p(x,a)=\psi(x)$
for all $a\in\A$.

\demo{Definition {\rm 5}}
A subset of $A$ of $Y$ will be called  $P$-{\it invariant}
 if $P\chi_A=\chi_A$
(where $\chi_A$ stands for the characteristic function of $A$).
\enddemo

Ergodicity of a measure-preserving Markov operator 
is equivalent to the absence of nontrivial invariant subsets  
(see \cite{rsbl}).
Lemma \ref{erg} follows now from:
 
\specialnumber{2}\proclaim{Proposition}\label{ergod}
Suppose that $\Pi$ is strictly irreducible{\rm .}
Then $A\subset Y$ is $P$\/{\rm -}\/invariant  if and only if 
$\chi_A$ does not depend on $\A$ and is $F_m$-invariant{\rm .}
\endproclaim

{\it Proof}. If  $p_{kl}>0$ then $\chi_A(T_kx,l)=\chi_A(x,k)$.
If $(\Pi^T\Pi)_{ij}>0$ then there exists $k\in\A$ such that 
$p_{ki}>0$, $p_{kj}>0$. Therefore, 
$\chi_A(x,k)=\chi_A(T_kx,i)=\chi_A(T_kx,j)$, which implies 
$\chi_A(x,i)=\chi_A(x,j)$ and proves that $\chi_A$ does not depend on $\A$.
The equality $\chi_A(T_ix,j)=\chi_A(x,i)$, true when  $p_{ij}>0$,
and the irreducibility of $\Pi$ imply group-invariance of $\chi_A$.

\specialnumber{3}\proclaim{Lemma}
\label{squerg}
Suppose that $\Pi$ is strictly irreducible and 
 $F_m^2$ acts ergodically on $(X,\nu)${\rm .}
Then the operator $P^2$ is ergodic{\rm .}
\endproclaim

By Lemma \ref{erg}, $P$ is ergodic. If $P^2$ is not ergodic, 
then there exists a nonconstant function $\psi\in \ly$ such that 
$P\psi=-\psi$. Arguing in the same way as in 
Proposition \ref{ergod}, 
we obtain that $\psi$ does not depend on $\A$, in other words,
there exists $\p\in\lx$ such that $\psi(x,a)=\p(x)$. 
The relation $P\psi=-\psi$ implies  $T_i\p=-\p$ for all $i\in\A$, whence
$T_g\p=\p$ for all $g\in F_m^2$, and the Lemma is proved.

\demo{{R}emark} The Kakutani-Hopf ergodic theorem
for Markov operators immediately implies that if the action of 
$F_m$ on $(X,\nu)$ is ergodic then for any $\p\in\lx$,
$$
\frac1N\sum\limits_{n=0}^{N-1}\sn\p\to\int\limits_X \p d\nu
$$ 
both $\nu$-almost everywhere and in $\lx$ as $N\to\infty$ 
(see \cite{gri2}, \cite{buf}).
\vglue6pt

The operator adjoint to $P$ is given by
\begin{equation}
\label{pstar}
P^*\p(x,i)= \sum_{j\in\A}\frac{p_jp_{ji}}{p_i}\p(T_{-j}x,j).
\end{equation}

Consider a unitary operator $U$ given by 
\begin{equation}
\label{uu}
U\p(x,i)=\p(T_ix,-i).
\end{equation}
Clearly, $U^2={\rm Id}$.

\specialnumber{3}\proclaim{Proposition}\label{conj}
Suppose  the matrix $\Pi$ satisfies the symmetry condition 
{\rm (\ref{sym}).} Then $P=UP^*U${\rm . }
\endproclaim

Indeed,  using (\ref{sym}), we can write 
$$
P^*\p(x,i)= \sum_{j\in\A}p_{-i,-j}\p(T_{-j}x,j)=UPU\p(x,i).
$$

\section{Uniform spherical averages}
\label{secunif}
\advance\eqcount by 8

In this section, we illustrate the method of the proof of 
Theorem \ref{sphfr}, by deducing the Nevo-Stein theorem 
from Rota's ``Alternierende Verfahren'' theorem \cite{rota} applied
to the Markov operator (\ref{operator}).

Consider uniform spherical averages (\ref{unif}). They are 
a particular case of the averages $\sn$ for $\Pi$ 
defined by  $p_{ij}=1/(2m-1)$ for $i+j\neq 0$ and 
$p_{ij}=0$ for $i+j=0$.

For $\Pi$ thus defined, the Markov operator (\ref{operator}) 
takes the form 
\begin{equation}
\label{punif}
P\p(x,i)=\frac 1{2m-1}\sum_{j:i+j\neq 0}\p(T_ix,j)
\end{equation}
and its adjoint is given by 
$$
P^*\p(x,i)=\frac 1{2m-1}\sum_{j:i+j\neq 0}\p(T_{-j}x,j).
$$

\specialnumber{4}\proclaim{Lemma}
\label{ident}
For $P$ given by {\rm (\ref{punif})} and $U$ given by {\rm (\ref{uu}),} 
$$
P^*P=\frac{2m-2}{2m-1}UP+\frac1{2m-1}{\rm Id}.
$$
\endproclaim

{\it Proof}.  We have 
$$
UP\p(x,i)=P\p(T_ix,-i)=\frac1{2m-1}\sum_{k:k\neq i}\p(x,k)
$$
and 
\begin{eqnarray*}
P^*P\p(x,i)&=&
\frac1{2m-1}\sum_{j:i+j\neq 0}P\p(T_{-j}x,j)\\
&=&
\frac 1{(2m-1)^2}
\sum_{j:i+j\neq 0}\sum_{k:j+k\neq 0}\p(x,k)\\
&=
&\frac{2m-2}{(2m-1)^2}\sum_{k:k\neq i}\p(x,k)+\frac1{2m-1}\p(x,i)\\
&=&
(\frac{2m-2}{2m-1}UP+\frac1{2m-1}{\rm Id})\p(x,i).
\end{eqnarray*}

From Lemma \ref{ident} and Proposition \ref{conj}, 
by induction, we obtain 
\begin{eqnarray}
(P^*)^nP^n&=&\frac{2m-2}{2m-1}UP^{2n-1}+\frac1{2m-1}(P^*)^{n-1}P^{n-1} 
\label{pncon}
\\
\noalign{\noindent or}
P^{2n-1}&=&\frac{2m-1}{2m-2}U(P^*)^nP^n-\frac1{2m-2}U(P^*)^{n-1}P^{n-1}.\nonumber
\end{eqnarray}

The Nevo-Stein theorem easily follows now from 
the Alternierende Verfahren theorem of Gian-Carlo Rota \cite{rota}:

\proclaimtitle{Rota \cite{rota}}
\specialnumber{4}\proclaim{Theorem}
\label{buro}\hskip-9pt 
Let $(Z,\mu)$ be a probability space{\rm .} Let 
$Q$ be a meas\-ure\/{\rm -}\/preserving  Markov operator on 
$L_1(Z,\mu)${\rm .}
Then for any $\p\in L\log L(Z,\mu)$ the sequence
$(Q^*)^nQ^n\p$ converges $\mu$\/{\rm -}\/almost everywhere
and in $L_1$ as $n\to\infty${\rm .}
\endproclaim

Theorem \ref{buro} generalizes Stein's theorem \cite{stein61}
on convergence of powers of self-adjoint operators and easily follows 
from the Martingale convergence
theorem; we recall its proof in Section \ref{markop}. Ornstein's  
counterexample \cite{ornstein} shows that 
neither Stein's nor Rota's theorem holds for $\p\in L_1$.

The equation (\ref{pncon})
and Theorem \ref{buro} yield the convergence of $P^{2n}\p$ for 
$\p\in L\log L(Y, \eta)$.
Lemma \ref{squerg} implies $F_m^2$-invariance of the limit. 
The Nevo-Stein theorem is proved.

\section{Proof of Lemma \ref{ppow}}
\label{secgencase}
\advance\eqcount by 10

\specialnumber5{}\proclaim{Lemma}
\label{maxin}
Suppose $\Pi$ generates the free group and satisfies the symmetry 
condition
{\rm (\ref{sym}).}
Then there exists a positive constant $c$ depending only on $\Pi$
such that 
for any nonnegative $\p\in\ly$ and any $n>0${\rm ,}
\begin{equation}
\label{comp}
(P^*)^nP^n\p\geq cUP^{2n-1}\p.
\end{equation}
\endproclaim

{\it Proof}.
We first prove the statement for $n=1$:
\begin{equation}
\label{compone}
P^*P\p\geq cUP\p .
\end{equation}
Now,
$$
P^*P\p(x,i)=\sum\limits_{j,k\in\A}\frac{p_jp_{ji}}{p_i}p_{jk}\p(x,k).
$$
If $\Pi$ generates the free group, then for any $i,k\in\A$ we have
$\sum_{j\in\A}\frac{p_jp_{ji}}{p_i}p_{jk}>0$.
Since 
$$
UP(x,i)=\sum\limits_k p_{-i,k}\p(x,k),
$$
(\ref{compone}) is proved; in view of Proposition \ref{conj},
(\ref{comp}) follows by induction, and the lemma is proved.
\vglue6pt

Now we  prove $L_1$-convergence of the powers $P^{n}$.
The following proposition is well known (see, for example, \cite{kaim}).

\specialnumber{4}\proclaim{Proposition}\label{tailconv}
Let $Q$ be a measure\/{\rm -}\/preserving Markov operator on a probability space
$(Z,\mu)${\rm .}
Then the tail sigma\/{\rm -}\/algebra of $Q$ is trivial if and only if 
for any $\p\in L_1(Z,\mu)${\rm ,}  $(Q^*)^n\p\to\int_Z \p d\mu$ 
in $L_1(Z,\mu)$ as $n\to\infty${\rm .} 
\endproclaim

Since $P=UP^*U$, triviality of the tail sigma-algebra of $P$ is equivalent 
to the triviality of that of $P^*$.
To establish this triviality, we shall use the following version of the 
$0$-$2$ law for Markov operators. 

\specialnumber{6}\proclaim{Lemma}
\label{trivtail}
Let $Q$ be an arbitrary measure\/{\rm -}\/preserving Markov operator on a 
probability space $(Z,\mu)${\rm .}

If the tail sigma\/{\rm -}\/algebra of $Q$ is trivial 
then for any $\p,\psi\in L_2(Z,\mu)$
$$\int\limits_Z (Q^*)^n\p\cdot (Q^*)^n\psi d\mu 
\to \int\limits_Z \p d\mu \int\limits_Z \psi d\mu$$
as $n\to\infty${\rm .}

If the tail sigma\/{\rm -}\/algebra of $Q$ is nontrivial then 
for any $\varepsilon>0$ there exist
positive functions $\p, \psi\in L_{\infty}(Z,\mu)$  
of integral $1$ such 
that 
$$\limsup\limits_{n\to\infty}\int(Q^*)^n\p\cdot(Q^*)^n\psi d\mu <\varepsilon.$$
\endproclaim

\vglue-24pt The proof of Lemma \ref{trivtail} closely models 
Vadim A. Kaimanovich's proof 
of the $0$-$2$ law \cite{kaim} and will be given in Section \ref{markop}.

\specialnumber{7}\proclaim{Lemma}
\label{mix}
Under assumptions of Lemma {\rm \ref{ppow},}
for any $\p, \psi\in L_2(Y, \eta)$,{\rm }
$$
\int_YP^n\p\cdot \psi d\eta \to \int_Y\p d\eta \int_Y\psi d\eta .
$$
\endproclaim

This follows from the $K$-property for the operator $P$, which we
prove in Section \ref{kprop} (Lemma \ref{KP}). 

Lemma \ref{trivtail}, Lemma \ref{mix}, and the inequality (\ref{comp}) 
easily imply triviality of the tail sigma-algebra of $P$.

Indeed, for any positive 
$\p,\psi\in L_{\infty}(Y,\eta)$, we have 
$$
\int_Y P^n\p\cdot P^n\psi d\eta=
\int_Y (P^*)^nP^n\p\cdot \psi d\eta 
\geq c \int_Y UP^{2n-1}\p\cdot \psi d\eta\to 
c\int_Y\p d\eta\int_Y\psi d\eta 
$$
as $n\to\infty$.
In view of Lemma \ref{trivtail}, this relation implies that 
$P^*$ (and hence also $P$, since $P=UP^*U$) 
has trivial tail sigma-algebra.

Proposition \ref{tailconv}
yields that for any $\p\in L_1(Y,\eta)$,
$$
P^n\p\to\int_Y\p d\eta
$$
in $L_1$ as $n\to\infty$.

Now we establish pointwise convergence of $P^n\p$ for $\p\in L\log L(Z,\mu)$.

Recall that if $(Z,\mu)$ is an arbitrary probability space 
then the {\it Orlicz norm} (see \cite{zygmund}, \cite{stein})
on the space $L\log L(Z,\mu)$ 
can be introduced, for example, 
by putting

$$||\p||_{L\log L}=\inf \{c: \int_Z \frac{|\p|}{c}\cdot 
\log (\frac{|\p|}{c}+2)d\mu\leq 1\}. 
$$

\specialnumber{8}\proclaim{Lemma}
\label{maxineq}
Let $(Z,\mu)$ be a probability space 
and let  $Q$ be a measure\/{\rm -}\/preserving Markov operator on 
$L_1(Z,\mu)${\rm .}

For any $p>1$ there exists a constant $A_p>0$
such that for any $\p\in L_p(Z,\mu)$ we have 
$$
||\sup\limits_n(Q^*)^nQ^n\p||_{L_p}\leq A_p||\p||_{L_p}.
$$
There exists a constant $A>0$ such that  
for any $\p\in L \log L(Z,\mu)${\rm ,}
$$
||\sup\limits_n(Q^*)^nQ^n\p||_{L_1}\leq A||\p||_{L\log L}.
$$
\endproclaim 

Lemma \ref{maxineq} will be proved in Section \ref{markop}.

Lemmas \ref{maxin}, \ref{maxineq} yield:

\specialnumber{9}\proclaim{Lemma}
\label{supin}
Let $p>1${\rm .}  Then there exists a constant $p>1$ such that for any 
 $\p\in L_p(Y,\eta)${\rm ,} 
$$
||\sup\limits_n P^{2n}\p ||_{L_p}\leq A_p||\p||_{L_p}.
$$
There exists a constant $A>0$ such that for any
 $\p\in L\log L(y, \eta)${\rm ,}
\vglue2pt\hfill ${\displaystyle
||\sup\limits_n P^{2n}\p ||_{L_1}\leq A||\p||_{L\log L}.
}$\hfill
\endproclaim

\specialnumber{5}\proclaim{Proposition}\label{tailltwo}
Let $Q$ be a measure\/{\rm -}\/preserving Markov operator on a probability space
$(Z,\mu)${\rm .}
If the tail sigma\/{\rm -}\/algebra of $Q^*$ is trivial then
for any $\p\in L_2(Z,\mu)$ we have  $Q^n\p\to\int \p$ 
in $L_2$ as $n\to\infty${\rm .} 
\endproclaim
 
The proof is given in Section \ref{markop}.
\vglue9pt

Now let $\p\in L_2(Y,\eta)$, $\int_Y\p d\eta=0$. 
Then $||P^n\p||_{L_2}\to 0$ as $n\to\infty$ by Proposition \ref{tailltwo}.
By  Lemma \ref{supin}, for any 
positive integer $k$, we have  
$$
||\sup\limits_n P^{2n+2k}\p ||_{L_2}\leq A_2||P^{2k}\p||_{L_2},
$$
and the right part of the inequality tends to $0$ as $k\to\infty$.
This implies pointwise convergence of $P^{2n}\p$ for $\p\in L_2(Y, \eta)$, 
and, since we have $L_1$-convergence for the 
whole sequence $P^n\p$, we also have pointwise convergence for $P^n\p$ with 
$\p\in L_2(Y, \eta)$.

Since $L_2$ is dense in $L\log L$, pointwise convergence 
of $P^n\p$ for $\p\in L_2$ and the $L\log L$-maximal inequality 
of Lemma \ref{supin} yield pointwise convergence of  $P^n\p$ 
for any $\p\in L\log L$.
\vglue6pt

To complete the proof of Lemma \ref{ppow} and
Theorem \ref{sphfr}, 
it only remains to prove 
Lemmas  \ref{trivtail},\ref{mix}, \ref{maxineq} and
Proposition \ref{tailltwo}. We do so 
in the following two sections.

\section{Proofs of Lemmas \ref{trivtail},   \ref{maxineq} and of 
Proposition \ref{tailltwo}}
\label{markop}

Let $(Z,\mu)$ be a probability space and 
let $Q$ be an arbitrary measure-preserving Markov operator on $L_1(Z,\mu)$. 
Let 
$$
\zz=\{\z=(z_n), n\in {\Bbb Z}, z_n\in Z\}
$$
be the space of bi-infinite sequences of elements of $Z$
and let $\qz$ be the Markov measure on $\zz$ 
corresponding to the operator $Q$ and the stationary distribution $\mu$.
Let $\sigma_Q$ be the shift on $(\zz, \qz)$ given by 
$(\sigma_Q(\z))_n=(\z)_{n+1}$; clearly, $\sigma_Q$ 
preserves the measure $\qz$.

For any $k,m\in \{-\infty\}\cup{\Bbb Z}\cup \{\infty\}$, $k\leq m$, 
denote by ${\cal F}_k^m$ the sigma-algebra on $\zz$
generated by the random variables
$z_l$, $k\leq l\leq m$. In particular, $\F_k$ is the sigma-algebra generated
by $z_k$. We shall sometimes write $F_{\geq k}$ for $\F_k^{\infty}$
and $\F_{\leq k}$ for $\F_{-\infty}^k$.

If $\p\in L_1(Z,\mu)$ and
 $\Phi\in L_1(\zz,\qz)$ is given by $\Phi(\z)=\p(z_0)$, then 
$E(\Phi(\z)|{\cal F}_{-\infty}^{-n})=Q^n\p(z_{-n})$, and
$E(E(\Phi(\z)|{\cal F}_{-\infty}^{-n})|{\cal F}_0)=(Q^*)^nQ^n\p(z_{0})$.

Rota's theorem (Theorem \ref{buro}) and Lemma \ref{maxineq}
immediately follow now from the inverted Martingale dominated 
convergence theorem and the corresponding maximal inequalities 
(see \cite[Chap.\ IV, Props.\ 2-8, 2-10]{neveu}).
This argument implies, moreover, the following:
 
\specialnumber{6}\proclaim{Proposition}\label{rotid}
Suppose that the tail sigma\/{\rm -}\/algebra of $Q^*$ is trivial{\rm .} Then
for all $\p\in L\log L(Z,\mu)${\rm ,}
$\lim\limits_{n\to\infty}(Q^*)^nQ^n\p=\int\p d\mu${\rm .} 
\endproclaim

Proposition \ref{rotid} implies Proposition \ref{tailltwo}, because 
if the tail sigma-algebra of a Markov operator is trivial then for any 
$\p\in L_2(Z,\mu)$ satisfying $\int_Z\p d\mu=0$, we have 
$$
\int_Z (Q^n\p)^2d\mu=\int_Z (Q^*)^nQ^n\p\cdot\p d\mu\to 0
$$
as $n\to\infty$, by Proposition \ref{rotid}. 
\vglue4pt

Now we prove Lemma \ref{trivtail}.
The proof closely models Kaimanovich's proof of the 
0-2 law for Markov operators \cite{kaim}.

The first part of the lemma is a corollary of 
Proposition \ref{tailltwo}.
To prove the second part,
let $\F_{\infty}$ be the tail sigma-algebra of $Q$, that is,  
$\F_{\infty}=\wedge_{k>0} \F_{\geq k}$, and
assume there exists  $A\in \F_{\infty}$ such that $0<\qz(A)<1$.
Set 
$$
\Phi(\z)=\chi_A(\z)/\qz(A),\ \Psi(\z)=\chi_{(\zz\setminus A)}(\z)/
\qz(\zz\setminus A).
$$ 
Then $\Phi, \Psi$ are positive, bounded, tail-measurable, 
$\int \Phi d\qz=\int \Psi d\qz=1$, 
$\Phi\cdot\Psi=0$. 
Let $M$ be a constant such that $M>\Phi$, $M>\Psi$.
Set $\p_k(z_k)=E(\Phi(\z)|\F_{\leq k})$, 
$\psi_k(z_k)=E(\Psi(\z)|\F_{\leq k})$.
Clearly, $\p_k, \psi_k$ are positive and bounded from above by $M$. By the 
Martingale convergence theorem, $\p_k(z_k)\to \Phi(\z)$, 
$\psi_k(z_k)\to\Psi(\z)$ both $\qz$-almost everywhere  
and in  $L_1(\zz, \qz)$ as $k\to\infty$.

Since $\Phi, \Psi$ are $\F_{\infty}$-measurable, we have 
$$E(\p_k(z_k)|\F_{\infty})\to \Phi,\ E(\psi_k(z_k)|\F_{\infty})\to \Psi$$
both $\qz$- almost everywhere and in  
$L_1(\zz, \qz)$ as $k\to\infty$.

Choose $k$ in such a way that 
$$
\int_{\zz}E(\p_k(z_k)|\F_{\infty})E(\psi_k(z_k)|\F_{\infty})d\qz < \varepsilon.
$$
Clearly, $$E(\p_k(z_k)|\F_{\geq n+k})=(Q^*)^n\p(z_{n+k}),
\ E(\psi_k(z_k)|\F_{\geq n+k})=(Q^*)^n\psi(z_{n+k}).$$
Therefore, as $n\to\infty$,
\begin{eqnarray*}
&&\hskip-18pt \int_Z (Q^*)^n\p_k(z)\cdot (Q^*)^n\psi_k(z) d\mu=
 \int_{\zz}E(\p_k(z_k)|\F_{\geq n+k})E(\psi_k(z_k)|\F_{\geq n+k})d\qz
\\&&\hskip.5in \to 
\int_{\zz}E(\p_k(z_k)|\F_{\infty})E(\psi_k(z_k)|\F_{\infty})d\qz < \varepsilon,
\\
\noalign{\noindent and Lemma \ref{trivtail} is proved.}
\end{eqnarray*}

\vglue-36pt
 \section{$K$-property and the proof of Lemma \ref{mix}}
\label{kprop}
\advance\eqcount by 12

Let $\yz$ be the space of biinfinite sequences  of elements of $Y$:
$$
\yz=\{\y: \y=(y_n), n\in{\Bbb Z}, y_n\in Y\}.
$$

Let $\pz$ be the measure corresponding to the operator $P$ and the 
stationary distribution $\eta$, and let 
 $\sip$ be the shift on $(\yz, \pz)$.  In order to prove Lemma~7, it suffices to show that $\sigma_P$ is mixing.  To do so, we
establish the following

\specialnumber{10}\proclaim{Lemma}
\label{KP}
Assume $F_m^2$ acts ergodically on $(X,\nu)$  and
assume the matrix $\Pi$ is strictly irreducible{\rm .} 
Then the system $(\yz, \pz, \sip)$ has $K$\/{\rm -}\/property{\rm .}
\endproclaim

The proof is based on the Rohlin-Sinai theorem \cite{rosin}.
First, we give another realization of $\sip$.

Let $\Sia$ be the space of bi-infinite sequences of symbols of $\A$:
$$
\Sia=\{\omega:\ \omega=(\omega_n), n\in{\Bbb Z},\omega_n\in\A\}.
$$

Let $\sia:\Sia\to\Sia$ be the shift on $\Sia$. 
Let $\mu_{\Pi}$ be the $\sia$-invariant
Markov measure on 
$\Sia$ corresponding to the matrix $\Pi$ and its stationary distribution~$p$.
Consider the map $T:\Sia\times X\to\Sia\times X$ given by the formula
\begin{equation}
\label{skewpr}
T(\omega,x)=(\sia\omega,T_{\omega_0}x).
\end{equation}
 Clearly, the map $T$ preserves the measure $\mu_{\Pi}\times\nu$.

\specialnumber{11}\proclaim{Lemma}
\label{iso}
The systems $(\Sia\times X, \mu_{\Pi}\times \nu, T)$ and
$(\yz, \pz, \sip)$ are isomorphic{\rm .}
\endproclaim

{\it Proof}.
Let $\y\in\yz$. Then  $\y=(y_n)$, where 
$y_n\in Y$; that is, $y_n=(i_n, x_n)$, $i_n\in\A$, $x_n\in X$.
Set  $\omega(\y)=(i_n)$, $n\in{\Bbb Z}$ 
and $x(\y)=x_0$.
The map $F:\yz\to\Sia\times X$ given by 
$F(\y)= (\omega(\y), x(\y))$ produces 
the desired isomorphism ($F$ is invertible because
for $\pz$ -- almost all $\y\in\yz$, we have
$x_1=T_{i_0}x_0$, $x_2=T_{i_1}x_1$, $x_{-1}=T_{-i_{-1}}x_0$, etc.)

Now we establish the  $K$-property for the system 
$(\Sia\times X, \mu_{\Pi}\times \nu, T)$.
The proof  follows the method of Oseledets
\cite{oseled}.

As in last section, write $\F_k^m(\yz)$ for the $\sigma$-algebra
in $\yz$ generated by the random variables $y_l$, $k\leq l\leq m$;
write $\F_k^m(\Sigma_{\A})$ for the $\sigma$-algebra in $\sigca$ generated
by the random variables $\omega_l$, $k\leq l\leq m$; write 
$\F_{\geq k}$  instead of $\F_k^{\infty}$, $\F_{\leq k}$ instead of
$\F_{-\infty}^{k}$, and $\F_k$ instead of $\F_k^k$; finally, denote by
$\B(X)$ the $\sigma$-algebra of all $\nu$-measurable subsets  of $X$, by 
$\B(\sigca\times X)$ the $\sigma$-algebra of all
$\mu_{\Pi}\times\nu$-measurable subsets of $\sigca\times X$.
 
Let $\pi(T)$ be the Pinsker $\sigma$-algebra of $T$. 
We shall use the Rohlin-Sinai theorem \cite{rosin} 
to prove the triviality of
$\pi(T)$, and therefore the $K$-property\break of $T$.

Consider the $\sigma$-algebra $\G_+$= $\F_{\geq 0}(\sigca)\times \B(X)$
(the future of our Markov process). Clearly,
$T\G_+\supset\G_+$, and $\vee_{k\in{\Bbb Z}}T^k\G_+=\B(\siax)$; by the
Rohlin-Sinai theorem [18], $\G_+\supset \pi(T)$. Let
$\G_-=\F_{\leq 0}(\sigca)\times \B(X)$ (the past of our Markov process).
 Clearly, 
$T^{-1}\G_-\supset \G_-$, and $\vee_{k\in{\Bbb Z}}T^k\G_-=\B(\siax)$.
By the Rohlin-Sinai theorem, 
$\G_-\supset\pi(T^{-1})=\pi(T)$. We have, therefore, $\pi(T)\subset
\G_+\wedge\G_-$. 

It is easy to check that $(\mu_{\Pi}\times\nu)-$ almost surely we have
\begin{equation}
\label{futpast}
(\F_{\geq 0}(\sigca)\times \B(X))\wedge (\F_{\leq 0}(\sigca)\times
\B(X))=\F_0(\sigca)\times \B(X).
\end{equation}

For $k\in{\Bbb Z}$, let $\G_k=\F_k(\sigca)\times \B(X)$ 
(the moment $k$ of our Markov process). 
By (\ref{futpast}), $\pi(T)\subset \G_0$. Since $T\pi(T)=\pi(T)$
and $T^k\G_0=\G_k$, 
$$\pi(T)\subset \wedge_{k\in{\Bbb Z}}\G_k.$$ 

Now let $\p:\siax\to\R$ be $\pi(T)$-measurable. 
Then for any $k\in{\Bbb Z}$ there exists 
$\psi: Y\to\R$ such that 
$\p(\omega,x)=\psi_k(\omega_k,x)$. 

Since for all $k\in{\Bbb Z}$ we have  $E(E(\p|\G_k)|\G_0)=\p$, we obtain 
\begin{equation}
\label{tofro}
(P^*)^kP^k\p_0=P^k(P^*)^k\p_0=\p_0\ \hbox{for all}\ k\in{\Bbb N}.
\end{equation}

To prove the  triviality of $\pi(T)$, it remains to prove that a function
$\p_0: Y\to\R$, satisfying (\ref{tofro}), is a constant.

\specialnumber{7}\proclaim{Proposition}
Suppose $\Pi$ is strictly irreducible{\rm ,} $\p\in\ly${\rm .} 
Then a set $A$ is $P^*P$\/{\rm -}\/invariant  if and only if 
$\chi_A$ does not depend on $\A${\rm .}
\endproclaim

Indeed, 
$$
P^*P\chi_A(x,i)=\sum_{k,l} \frac{p_kp_{ki}}{p_i}p_{kl}\chi_A(x,l)
$$
and, for  $i,l$ fixed, we have  
$\sum_k\frac{p_kp_{ki}}{p_i}p_{kl}>0$
if and only if $(\Pi^T\Pi)_{il}>0$, which implies the proposition.

\specialnumber{12}\proclaim{Lemma}
\label{pdva}
Suppose the matrix $\Pi$ is strictly irreducible{\rm .} 
Suppose a set $A\subset Y$ is both $P^*P$ and $(P^*)^2P^2$ 
invariant{\rm .} Then $\chi_A$ does not depend on $\A$ 
and is $F_m^2$\/{\rm -}\/invariant{\rm .} 
\endproclaim

By the previous proposition, $\chi_A$ does not depend on $\A$. 
Write
$$\chi_A(x,i)=(P^*)^2P^2\chi_A(x,i)=
\sum\frac{p_jp_{ji}}{p_i} \frac{p_kp_{kj}}{p_j}p_{kl}p_{lm}\chi_A(T_lT_{-j}x,m).
$$
We have, then,
$\chi_A(x)=\chi_A(T_lT_{-j}x)$ for all $j,l$
such that $(\Pi^T\Pi)_{jl}>0$. Since the matrix
$\Pi$ is strictly irreducible, 
the claim is proved.

Lemma \ref{KP} is proved and it implies, in particular, that 
$\sigma_P$ is mixing, which yields Lemma \ref{mix}.

The proof of Theorem \ref{sphfr} is complete.

\demo{{R}emark {\rm 1}} Let $(Z,\mu)$ be a Lebesgue 
probability space, $Q$ a measure-preserving Markov operator on $L_1(Z,\mu)$,
$(\zz, \qz)$ the space of trajectories of $Q$, $\sigma_Q$ 
the corresponding shift, and $\pi(\sigma_Q)$ the Pinsker sigma-algebra
of $\sigma_Q$. Then we have:

\specialnumber{8}\proclaim{Proposition}\label{pins}
 Note that $\pi(\sigma_Q)\subset \F_0(\zz)${\rm .}
If $C\subset Z$ and the set $\{\z: z_0\in C\}\in \pi(\sigma_Q)$ then
$\chi_C=Q^k(Q^*)^k\chi_C=
(Q^*)^kQ^k\chi_C$ for any $k\in{\Bbb N}${\rm .} 
\endproclaim

The proof is the same as that of Lemma \ref{KP}:
first, the Rohlin-Sinai theorem gives that 
$\pi(\sigma_Q)\subset \F_{\geq 0}(\zz)\cap \F_{\leq 0}(\zz)=\F_0(\zz)$, 
then the $\sigma_Q$-invariance of $\pi(\sigma_Q)$ implies that 
$\pi(\sigma_Q)\subset \wedge_{k\in{\Bbb Z}}\F_k(\zz)$, which implies 
the proposition. 
\enddemo

 {\it Remark} 2.  Let $\mu$ be an arbitrary Borel probability
$\sigma_{\A}$-invariant measure on
$\sigca$. Clearly, the map $T$, defined by (\ref{skewpr}), 
preserves the measure $\mu\times\nu$. 

Let $B=(B_{ij})$, $i,j\in \A$, be a $0-1$ matrix, and 
let $\mu$ be a Gibbs measure (in the sense of Bowen \cite{bowen})
on the subshift of $\Sigma_{\A}$ given by the matrix $B$. 
\vglue6pt

Arguing in the same way as in the proof of Lemma \ref{KP}, 
we see that if $B$ is strictly irreducible and
the action of $F_m^2$ on $X$ is ergodic, then the  
system $(\siax, \mu\times\nu, T)$ has the $K$-property.

\section{Proof of Theorem \ref{sphir}}
\label{endof}
\specialnumber{13}\proclaim{Lemma}
\label{mair}
Suppose that $\Pi$ satisfies {\rm (\ref{sym})}
and that all entries of the matrix
$\Pi\Pi^T+(\Pi\Pi^T)^2+\cdots +(\Pi\Pi^T)^k$ are positive{\rm .} 
Then there exists 
a constant $c>0$ such that for any nonnegative 
$\p\in\ly${\rm ,}
$$
((P^*)^nP^n+(P^*)^nPP^*P^n+\dots+(P^*)^n(PP^*)^kP^n)\p
\geq cUP^{2n-1}\p  
$$
almost everywhere{\rm .}
\endproclaim

The proof is the same as that of Lemma \ref{maxin}.

\specialnumber{14}\proclaim{Lemma}
\label{qqst}
Let $(Z,\mu)$ be a probability space{\rm ,} let  
$Q$ be a measure\/{\rm -}\/preserving Markov operator on 
$L_1(Z,\mu)${\rm ,} let $p>1$ and let $k$ be a positive integer{\rm .}
Then for any $\p\in L_p(X,\nu)$ the sequence
$Q^n(Q^*Q)^k(Q^*)^n\p$ converges $\mu$\/{\rm -}\/almost everywhere
and in $L_p$ as $n\to\infty${\rm .}

Moreover{\rm ,} 
$$
||\sup\limits_n(Q^*)^n(Q^*Q)^kQ^n\p||_{L_p}\leq A_p||\p||_{L_p}.
$$
\endproclaim

Let $\p\in L_1(Z,\mu)$ and
define $\Phi\in L_1(\zz,\qz)$ by $\Phi(\z)=\p(z_0)$.
Set 
$$
\Phi_n^0(\z)=E(\Phi(\z)|{\cal F}_{n}^{\infty})=
(Q^*)^n\p(z_{n})
$$
and for $i\geq 1$ let 
$$
\Phi_n^i(\z)=E(E(\Phi_n^{i-1}(\z)|\F_{n-1})|\F_n).
$$
Clearly, 
$$\Phi_n^i(\z)=(Q^*Q)^i(Q^*)^n\p(z_{n}),$$ 
and $$E(\Phi_n^i(\z)|\F_0)=Q^n(Q^*Q)^i(Q^*)^n\p(z_{n}).$$

The statement of the proposition follows now from 
Rota's theorem (Theorem \ref{buro}) and the $L_p$ maximal inequality
for martingales (see \cite[Prop.\ IV-2-8]{neveu})
by induction on $i$.

In a similar fashion,  Lemma \ref{trivtail} implies:
\specialnumber{15}\proclaim{Lemma} 
\label{trivqq}
Let $k$ be a nonnegative integer{\rm .}
If the tail sigma-algebra of $Q$ is trivial 
then for any $\p,\psi\in L_2(Z,\mu)$
$$\int Q^k(Q^*)^n\p\cdot  Q^k(Q^*)^n\psi d\mu \to \int \p d\mu \int \psi d\mu$$
as $n\to\infty${\rm .}

If the tail sigma\/{\rm -}\/algebra of $Q$ is nontrivial then 
for any $\varepsilon>0$ there exist
positive functions $\p, \psi\in L_{\infty}(Z,\mu)$  
of integral $1$ such 
that 
$$\limsup_{n\to\infty}\int Q^k(Q^*)^n\p\cdot Q^k(Q^*)^n\psi d\mu <\varepsilon.$$
\endproclaim

The rest of the proof goes the same way as that of Theorem \ref{sphfr}, 
with Lemma \ref{qqst} being used instead of 
Lemma \ref{maxineq} and Lemma \ref{trivqq} assuming the role of Lemma~\ref{trivtail}.

\section{A conjecture}

Theorem \ref{sphir} can be applied to obtain spherical 
convergence for actions of some classes of Markov groups 
(in the sense of Gromov \cite{gromov}). 

Let $\Gamma$ be a Markov group. Its elements can then be coded 
by admissible words in a topological Markov chain. 
Assume that the matrix $A$ of the chain is irreducible and let $\Pi$ 
be the matrix of the Parry measure 
(in other words, the measure of maximal entropy) 
corresponding to $A$.
 If $\Pi$ is strictly irreducible
and satisfies the symmetry condition (\ref{sym}), then Theorem \ref{sphir}
is applicable. The spherical averages $\sn$ for $\Pi$ thus chosen 
can easily be reduced to {\it uniform} spherical averages in $\Gamma$
(see \cite{rokhfin}). 
Theorem \ref{sphir} then yields convergence of uniform 
spherical averages for the group $\Gamma$. 
For example, this takes place for Vershik's locally 
finite groups \cite{vershik}.

Gromov \cite{gromov} proved that Gromov hyperbolic groups 
are Markov. If the coding satisfied the assumptions of Theorem \ref{sphir},
then Theorem \ref{sphir} would yield the following:

\specialnumber{1}\proclaim{{C}onjecture}
Let $\Gamma$ be a Gromov hyperbolic group{\rm ,} let 
$S$ be a symmetric set of generators{\rm ,} and denote by $\Gamma^2$
the subgroup generated by elements that have a geodesic representation of even length over the alphabet $S$.
Let $p>1${\rm .}
Suppose $\Gamma$ acts 
on a probability space $(X,\nu)$ by measure\/{\rm -}\/preserving transformations{\rm .}

Then for any $\p\in L_p(X,\nu)$ the sequence
$$
s_{2n}\p=\frac1{\#\{g:|g|_S=2n\}}\sum_{g: |g|_S=2n}T_g\p
$$
converges as $n\to\infty$ almost everywhere and in $L_p$ to 
a $\Gamma^2$\/{\rm -}\/invariant function{\rm .}
\endproclaim

Assuming exponential mixing, 
Fujiwara and Nevo \cite{fune} obtained a convergence theorem for 
Cesaro averages of the spherical averages for Gromov hyperbolic groups.

\demo{Acknowledgements} 
I am deeply grateful to 
Rostislav I.\ Grigorchuk, who introduced 
me to this subject, and to 
Amos Nevo, who  suggested to use\break  Rota's theorem. 
Vadim A.\ Kaimanovich made many important suggestions, 
both  on  content and on presentation; I am greatly indebted to him. 
I am\break grateful to Charles Fefferman, Boris M.\ Gurevich, Gregory A.\
Margulis,\break  Valeriy I.\ Oseledets, Yakov G.\ Sinai, Elias M.\ Stein, 
Jean-Paul Thouvenot, and Anatoly M.\ Vershik for useful discussions.
While I was working on the paper, I visited La Sapienza di Roma and 
KTH Stockholm in the
framework of ``Russian-Swedish Workshop on Dynamical Systems".
I am deeply grateful to these institutions for
their hospitality.
This research was partially supported by the CRDF under grant RM1-2086.
\enddemo

 \end{document}